\documentclass[a4paper,12pt]{article}
\usepackage{a4wide}
\usepackage{amsmath}
\usepackage{amssymb}
\usepackage{amsthm}
\usepackage{latexsym}
\usepackage{graphicx}
\usepackage[english]{babel}
\usepackage{makeidx}

\newtheorem{obs} [subsection]{Remark}
\newtheorem{exm} [subsection]{Example}

\newtheorem{prop}[subsection]{Proposition}

\newtheorem{teor}[subsection]{Theorem}
\newtheorem{lema}[subsection]{Lemma}
\newtheorem{cor} [subsection]{Corollary}
\newcommand{\Zng}{finitely generated $\mathbb Z^n$-graded $S$-module}
\newcommand{\me}{\mathbf m}
\def\sdepth{\operatorname{sdepth}}
\def\depth{\operatorname{depth}}
\begin{document}
\selectlanguage{english}
\frenchspacing

\large
\begin{center}
\textbf{Some remarks on the Stanley depth for multigraded modules.}

Mircea Cimpoea\c s
\end{center}
\normalsize

\begin{abstract}
We show that Stanley's conjecture holds for any multigraded $S$-module $M$ with $\sdepth(M)=0$, where $S=K[x_1,\ldots,x_n]$. Also, we give some bounds for the Stanley depth of the powers of the maximal irrelevant ideal in $S$.
\vspace{5 pt} 


\vspace{5 pt} \noindent \textbf{2000 Mathematics Subject
Classification:}Primary: 13H10, Secondary: 13P10.
\end{abstract}

\section*{Introduction}

Let $K$ be a field and $S=K[x_1,\ldots,x_n]$ the polynomial ring over $K$. Let $M$ be a \Zng. A \emph{Stanley decomposition} of $M$ is a direct sum $\mathcal D: M = \bigoplus_{i=1}^rm_i K[Z_i]$ as $K$-vector space, where $m_i\in M$, $Z_i\subset\{x_1,\ldots,x_n\}$ such that $m_i K[Z_i]$ is a free $K[Z_i]$-module. The latter condition is needed, since the module $M$ can have torsion. We define $\sdepth(\mathcal D)=min_{i=1}^r |Z_i|$ and $\sdepth(M)=max\{\sdepth(M)|\;\mathcal D$ is a Stanley decomposition of $M\}$. The number $\sdepth(M)$ is called the \emph{Stanley depth} of $M$. Herzog, Vladoiu and Zheng show in \cite{hvz} that this invariant can be computed in a finite number of steps if $M=I/J$, where $J\subset I\subset S$ are monomial ideals. A computer implementation of this algorithm, with some improvements, is given by Rinaldo in \cite{rin}.

Let $M$ be a \Zng. Stanley's conjecture says that $\sdepth(M)\geq \depth(M)$. The Stanley conjecture for $S/I$ was proved for $n\leq 5$ and in other special cases, but it remains open in the general case. See for instance, \cite{apel}, \cite{hsy}, \cite{jah}, \cite{pops}, \cite{popi} and \cite{pop}. Another interesting problem is to explicitly compute the sdepth. This is difficult, even in the case of monomial ideals! Some small progresses were made in \cite{asia}, \cite{hvz}, \cite{mir}, \cite{mir2} and \cite{shen}.

In the first section, we prove that the Stanley conjecture holds for modules with $\sdepth(M)=0$, see Theorem $1.4$. As a consequence, it follows that any torsion free module $M$ has $\sdepth(M)\geq 1$. In the second section, we give an upper bound for the Stanley depth of the powers of the maximal ideal $\me=(x_1,\ldots,x_n)\subset S$, see Theorem $2.2$. We conjecture that $\sdepth(\me^k)=\left\lceil \frac{n}{k+1}\right\rceil$, for any positive integer $k$.

\textbf{Aknowledgements}. The author would like to express his gratitude to the organizers of PRAGMATIC 2008, Catania, Italy and especially to Professor Jurgen Herzog.

\footnotetext[1]{This paper was supported by CNCSIS, ID-PCE, 51/2007}

\newpage
\section{Stanley's conjecture for modules with sdepth zero.}

Let $M$ be a \Zng. We use an idea of Herzog, in order to obtain a decomposition of $M$, similar to the Janet decomposition given in \cite{imran}. For any $j\geq 1$, we have a natural surjective map $\varphi_j: M \rightarrow x_n^j M$ given by the multiplication with $x_n^j$. Obviously, $\varphi_j(x_nM) \subset x_n^{j+1}M$ and therefore $\varphi_j$ induces a natural surjection $\bar{\varphi}_j: M/x_nM \rightarrow x_n^jM/ x_n^{j+1}M$. We write $L_j=Ker (\bar{\varphi}_j)$. 

Note that $L_j\subset L_{j+1}$ for any $j$, since we have a natural surjection $x_n^jM/ x_n^{j+1}M \rightarrow x_n^{j+1}M/ x_n^{j+2}M$ given by multiplication with $x_n$. As $M/x_nM$ is finitely generated, it follows that there exists a nonnegative integer $q$ such that $L_{q}=L_{q+1}=\cdots$ and moreover $x_n^jM/ x_n^{j+1}M \cong x_n^{j+1}M/ x_n^{j+2}M$ for any $j\geq q$. Now, we can prove the following Lemma.

\begin{lema}
Let $M$ be a \Zng \hspace{5pt} and $q$ such that $L_{q}=L_{q+1}=\cdots$. Then we have the following decomposition of $M$ , as $K$-vector space: \[ M \cong M/x_nM \oplus \cdots \oplus x_n^{q-1}M/ x_n^{q}M \oplus x_n^{q}M/ x_n^{q+1}M [x_n]. \]
\end{lema}

\begin{proof}
Note that, since $M$ is graded, $\bigcap x_n^j M = 0$. Therefore, we have \[ M = M/x_nM \oplus x_nM = M/x_nM \oplus x_nM/x_n^2M \oplus x_n^2M = \cdots = \bigoplus_{j\geq 0} x_n^jM/x_n^{j+1}M.\] Since $x_n^jM/ x_n^{j+1}M \cong x_n^{j+1}M/ x_n^{j+2}M$ for any $j\geq q$, the proof of Lemma is complete.
\end{proof}

Note that each factor $x_n^{j}M/ x_n^{j+1}M$ naturally carries the structure of a multigraded $S'$-module, where $S'=K[x_1,\ldots,x_{n-1}]$.
Also, if $M=S/I$, where $I\subset S$ is a monomial ideal, the above decomposition is exactly the Janet decomposition of $S/I$, with respect to the variable $x_n$.

\begin{lema}
Let $M$ be a \Zng. Then $\sdepth(M)=n$ if and only if $M$ is free.
\end{lema}

\begin{proof}
If $M$ is free, it follows that $M \cong \bigoplus_{i=1}^r S(-a_i)$, where $a_i\in\mathbb Z^n$ are some multidegrees. Therefore, $M$ has a basis $\{e_1,\ldots,e_n\}$ where $e_i$ correspond to $1\in S(-a_i)$. Therefore \linebreak $M = \bigoplus e_i S$ is a Stanley decomposition of $M$ and thus $\sdepth(M)=n$. Conversely, given a Stanley decomposition $M = \bigoplus e_i S$, it follows that $M \cong \bigoplus_{i=1}^r S(-a_i)$, where $deg(e_i)=a_i$.
\end{proof}

\begin{lema}
Let $M$ be a graded $K[x]$-module. Then, the following are equivalent:

(1) $M$ is free.

(2) $M$ is torsion free.

(3) $\depth(M)=1$.

(4) $\sdepth(M)=1$.
\end{lema}

\begin{proof}
The equivalences $(1)\Leftrightarrow (2) \Leftrightarrow (3)$ are well known. $(4) \Leftrightarrow (1)$ is the case $n=1$ of the previous Lemma.
\end{proof}

Let $\me = (x_1,\ldots,x_n) \subset S$ be the maximal irrelevant ideal. Let $M$ be a \Zng. We denote $sat(M)=(0:_M \me^{\infty}) = \bigcup_{k\geq 1} (0:_M \me^{k})$ the \emph{saturation} of $M$. It is well known, that $\depth(M)=0$ if and only if $\me\in Ass(M)$ if and only if $sat(M)\neq 0$. On the other hand, $sat(M/sat(M))=0$. Note that if $I\subset S$ is a monomial ideal, then $sat(S/I)=I^{sat}/I$, where $I^{sat}=(I:\me^{\infty})$ is the saturation of the ideal $I$. We prove the following generalization of \cite[Theorem 1.5]{mir2}.

\begin{teor}
Let $M$ be a \Zng. If $\sdepth(M)=0$ then $\depth(M)=0$. Conversely, if $\depth(M)=0$ and $dim_K(M_a)\leq 1$ for any $a\in\mathbb Z^n$, then $\sdepth(M)=0$.
\end{teor}

\begin{proof}
We use induction on $n$. If $n=1$, then we are done by Lemma $1.3$. Suppose $n>1$.
We consider the decomposition \[ (*)\;\;M \cong M/x_nM \oplus \cdots \oplus x_n^{q-1}M/ x_n^{q}M \oplus x_n^{q}M/ x_n^{q+1}M [x_n],\] given by Lemma $1.2$. We define $M_j:=x_n^{j}M/ x_n^{j+1}M$ for $j\in [q]$. 
Since $\sdepth(M)=0$, it follows that $\sdepth(M_j)=0$ for some $j<q$. We have $M_j=sat(M_j)\oplus M/sat(M_j)$, where $sat(M_j)$ is the saturation of $M_j$ as a $S'$-module. If there exists some nonzero element $m\in sat(M_j)$ such that $x_n^j m = 0$, it follows that $m\in sat(M)$ and thus $sat(M)\neq 0$. 

For the converse, we assume $depth(M)>0$. It follows that $x_n sat(M_j)\subset sat(M_{j+1})$ for any $j<q$. Since $sat(M_j/sat(M_j))=0$, by induction hypothesis, it follows that $\sdepth(M_j/sat(M_j))\geq 1$. Therefore, $(*)$ implies
\[ (**) M \cong \bigoplus_{j=0}^{q-1} M_j/sat(M_j) \oplus M_q/sat(M_q) [x_n] \oplus \bigoplus_{j=0}^{q-1} sat(M_j) \oplus sat(M_q)[x_n].\]
On the other hand, $\bigoplus_{j=0}^{q-1} sat(M_j) \oplus sat(M_q)[x_n] = \bigoplus_{j=0}^{q} \bigoplus_{\bar{m}\in sat(M_j)/sat(M_{j-1})} m K[x_n]$ since $dim_K(M_a)\leq 1$, and therefore, by $(**)$, we obtain a Stanley decomposition of $M$ with it's $\sdepth\geq 1$! %For the converse, note that any nonzero element $u\in Soc(M)$ is annihilated by any monomial $z\in S$.
\end{proof}

\begin{cor}
If $M$ is torsion free, then $\sdepth(M)\geq 1$.
\end{cor}

\begin{proof}
Obviously, since $M$ is torsion free, we have $\depth(M)\geq 1$.
\end{proof}

\begin{exm}(Dorin Popescu, \cite{pop})
The condition $dim_K(M_a)\leq 1$ is essential in the second part of Theorem $1.4$. Let $S=K[x_1,x_2]$ and consider the module $M:=(Se_1\oplus Se_2)/(x_1z,x_2z$, where $z=x_1e_2-x_2e_1$. $M$ is multigraded with $deg(e_1)=deg(x_1)=(1,0)$ and $deg(e_2)=deg(x_2)=(0,1)$. Note that $dim_K(M_a)=1$ for any $a\in \mathbb Z^2\setminus \{(1,1)\}$ and $dim_K(M_{(1,1)})=2$. Since $z\in Soc(M)$, it follows that $\depth(M)=0$. We have a Stanley decomposition of $M$,
\[M = \bar{e}_1 K[x_2] \oplus \bar{e}_1 x_1K[x_1] \oplus \bar{e}_2 K[x_1] \oplus \bar{e}_2 x_2 K[x_2] \oplus\bar{e_1}x_1x_2K[x_1,x_2], \]
where $\bar{e_1}, \bar{e_2}$ are the images of $e_1$ and $e_2$ in $M$. It follows that $\sdepth(M)\geq 1$ and thus $\sdepth(M)=1$, since $M$ is not free.
\end{exm}

\begin{obs}
Let $M$ be a torsion free \Zng. Then we have an inclusion $0 \rightarrow M \rightarrow F$, where $F$ is a free module with $rank(F)=rank(M)$. Let $Q:=F/M$. Is it true that $\sdepth(M)\geq \sdepth(Q)+1$? 
%Reformulation: Given $Q$ a \Zng with $t(Q)=Q$ let $\varphi:F \rightarrow Q$ be a surjection of graded $\mathbb Z^n$ modules, with $F$ free. If $M:=Ker(\varphi)$, then $\sdepth(M)\geq \sdepth(Q)+1$. 
In particular, if $I\subset S$ is a monomial ideal, is it true that $\sdepth(I)\geq \sdepth(S/I)+1$?

If this result were true, then by $\depth(M)=\depth(Q)+1$, if $Q$ satisfy Stanley's conjecture, then $M$ also satisfy Stanley's conjecture. Note that, in general we cannot expect that $\sdepth(M)=\sdepth(Q)+1$. Take for instance $M=\me=(x_1,\ldots,x_n)\subset S$ and $Q=k=S/\me$. It is known from \cite{hvz} and \cite{par} that $\sdepth(\me)=\left\lceil \frac{n}{2} \right\rceil$, but $\sdepth(k)=0$. It would be interesting to characterize those modules $M$ with $\sdepth(M)=\sdepth(Q)+1$. Or, at least, the monomials ideals $I\subset S$ with $\sdepth(I)=\sdepth(S/I)+1$.
\end{obs}

We end this section with the following example.

\begin{exm}
Let $M_i:=syz_i(K)$ the $i$-th syzygy module of $K$. It is known that $\depth(M_i)=i$ for all $0\leq i\leq n$. The problem of computing $\sdepth(M_i)$ is a chellenging problem. Obviously, $\sdepth(M_0)=\sdepth(K)=0$. On the other hand, $\sdepth(M_1)=\sdepth(\me)=\left\lceil \frac{n}{2} \right\rceil$. Also, $\sdepth(M_n)=\sdepth(S)=n$. We claim that $\sdepth(M_{n-1})=n-1$.

Indeed, $M_{n-1}=Coker(S \stackrel{\psi}{\longrightarrow} S^n)$, where we define $S^n = \bigoplus_{i=1}^n Se_i$ and $\psi(1):=x_1e_1+\cdots+x_ne_n$. Therefore, $M_{n-1}:=S\bar{e}_1 + \cdots + S\bar{e}_n$, where $\bar{e}_i$ are the class of $e_i$ in $M_{n-1}$ for all $i\in [n]$. Note that $\bar{e}_1,\ldots,\bar{e}_{n-1}$ are linearly independent in $M_{n-1}$, since the only relation in $M_{n-1}$ is $x_1\bar{e}_1+\cdots+x_{n-1}\bar{e}_n = - x_n \bar{e}_n$. It follows that, 
\[ M_{n-1} = S\bar{e}_1 \oplus \cdots \oplus S\bar{e}_{n-1} \oplus K[x_1,\ldots,x_{n-1}]\bar{e}_n, \]
and therefore $\sdepth(M_{n-1})\geq n-1$. On the other hand, $\sdepth(M_{n-1})\leq n-1$, since $M$ is not free. Thus $\sdepth(M_{n-1})=n-1$.
\end{exm}

\section{Bounds for the sdepth of powers of the maximal irrelevant ideal}

Let $\mathbf{m}=(x_1,\ldots,x_n)$ be the maximal irrelevant ideal of $S$. Let $k\geq 1$ be an integer. In this section, we will give some upper bounds for $\sdepth(\me^k)$. In order to do so, we consider the following poset, associated to $\me^k$, 
\[ P:= \{u\in \me^k \;monomial\;: \; u|x_1^kx_2^k\cdots x_n^k \}, \]
where $u\leq v$ if and only if $u|v$. For any $u\in P$, we denote $\rho(u)=|\{j: x_j^k|u\}|$. Note that, by \cite[Theorem 2.4]{hvz}, there exists a partition of $P=\bigoplus_{i=1}^r [u_i,v_i]$, i.e. a disjoint sum of intervals $[u_i,v_i]=\{u\in P:\; u_i|u$ and $u|v_i \}$, such that $\min_{i=1}^r\{\rho(v_i)\} = \sdepth(\me^k)$.

We write $P_d =\{ u\in P: deg(u)=d \}$, where $k\leq d\leq kn$, and
$\alpha_d:=|P_d|$. First, we want to compute the numbers $\alpha_d$.

\begin{lema}
We the above notations, we have:
\[ \alpha_d=\sum_{i\geq 0} (-1)^i \binom{n}{i} \binom{n+d-i(k+1)-1}{n-1}.\]
\end{lema}

\begin{proof}
We fix $d\geq k$. For any $j\in [n]$, we write $A_j:=\{u\in S:\;deg(u)=d$, $x_j^{k+1}|u \}$. Obviously, $P_d:=S_d \setminus (A_1\cup A_2\cup \cdots \cup A_n)$, where $S_d$ is the set of all monomials of degree $d$ in $S$. For any nonempty subset $I\subset [n]$, we write $A_I:=\bigcap_{i\in I}A_i$. By inclusion-exclusion principle,
\[ |A_1\cup \cdots \cup A_n| = \sum_{\emptyset\neq I\subset [n]}(-1)^{|I|-1}|A_I|. \]
Note that a monomial $u\in A_I$ can be written as $u=w\cdot \prod_{i\in I}x_i^{k+1}$. Therefore, $|A_I|=\binom{n+d-i(k+1)-1}{n-1}$. Now, one can easily get the required conclusion.
\end{proof}

\begin{teor}
Let $a\leq \left\lceil \frac{n}{2} \right\rceil$ be a positive integer. Then $\sdepth(\me^k)\leq \left\lceil \frac{n}{k+1} \right\rceil$.
%If $k\geq \left\lceil \frac{n-a}{a} \right\rceil$ then $\sdepth(\me^k)\leq a$. 
In particular, if $k\geq n-1$, then $\sdepth(\me^k)=1$.
\end{teor}

\begin{proof}
Let $a=\left\lceil \frac{n}{k+1} \right\rceil$ and assume, by contradiction, that $\sdepth(\me^k)\geq a+1$. Obviously, by Lemma $2.1$, $\alpha_{k}=\binom{n+k-1}{n-1}$ and $\alpha_{k+1}=\binom{n+k}{n-1}-n$. We consider a partition of $\mathcal P: P_{n,k}=\bigcup_{i=1}^r[x^{c_i},x^{d_i}]$ with $\sdepth(\mathcal D(\mathcal P)) = a+1$. Note that $\me^k$ is minimally generated by all the monomials of degree $k$ in $S$. We can assume that $S_k = \{ x^{c_i}|i=1,\ldots,N \}$, where $N=\binom{n+k-1}{n-1}$. We consider an interval $[x^{c_i},x^{d_i}]$. If $c_i=x_j^k$, then by $\rho(x^{d_i})\geq a+1$, it follows that in $[x^{c_i},x^{d_i}]$ are at least $a$ distinct monomials of degree $k+1$. If $c_i(j)<k$ for all $j\in [n]$, then, in $[x^{c_i},x^{d_i}]$ are at least $a+1$ distinct monomials of degree $k+1$.

We assume that $k\geq \left\lceil \frac{n-a}{a} \right\rceil$. Since $\mathcal P: P_{n,k}=\bigcup_{i=1}^r[x^{c_i},x^{d_i}]$ is a partition of $P_{n,k}$, by above considerations, it follows that $\alpha_{k+1}\geq na+ (\alpha_k-n)(a+1)$. Therefore, $\binom{n+k}{k-1} \geq (a+1) \binom{n+k-1}{n-1}$. This implies $n+k\geq (k+1)(a+1) \geq (k+1)(\frac{n}{k+1}+1) = n+k+1$, a contradiction.
\end{proof}

We conjecture that $\sdepth(\me^k)\leq \left\lceil \frac{n}{k+1} \right\rceil$. Using the computer, see \cite{rin}, one can prove that this conjecture is true for small $n$. Also, the conjecture is true for $k=1$, from \cite{hvz}, \cite{par}. We end this section with the following proposition.

\begin{prop}
Let $I\subset S$ be a monomial ideal. Then $\sdepth(\me^k I)=1$ for $k\gg 0$.
\end{prop}

\begin{proof}
We consider the $K$-algebra $A:=\bigoplus_{i\geq 0} \me^i I/\me^{i+1} I$ and denote $A_i$ the $i^{th}$ graded component of $A$. Note that $H(A,i):=dim_K(A_i)=|G(\me^i I)|$, where $G(\me^i I)$ is the set of minimal monomial generators of $\me^i I$. Since $A$ is a finitely generated $K$-algebra, it follows that the Hilbert function $H(A,i)$ is polynomial for $i\gg 0$. 

Therefore, $lim_{i\rightarrow \infty} H(A,i)/H(A,i+1) = 1$. Note that there are exactly $H(A,i+1)$ monomials of degree $i+1$ in $\me^i I$. Suppose $\sdepth(\me^i I)\geq 2$. As in the proof of Theorem $2.2$, it follows that $H(A,i+1) \geq 2(H(A,i)-n)+n$, which is false for $i\gg 0$, since it contradicts the fact that $\lim_{i\rightarrow \infty} H(A,i)/H(A,i+1) = 1$.
\end{proof}

 \vspace{2mm} \noindent {\footnotesize
\begin{minipage}[b]{15cm}
 Mircea Cimpoeas, Institute of Mathematics of the Romanian Academy, Bucharest, Romania\\
 E-mail: mircea.cimpoeas@imar.ro
\end{minipage}}

\begin{thebibliography}{99}%{breitestes Label}  
  \bibitem[1]{pops}Sarfraz Ahmad, Dorin Popescu "Sequentially Cohen-Macaulay monomial ideals of embedding dimension four", Bull. Math. Soc. Sc. Math. Roumanie 50(98), no.2 (2007), p.99-110.
  \bibitem[2]{imran}Imran Anwar "Janet's algorithm", Bull. Math.
                    Soc. Sc. Math. Roumanie 51(99), no.1 (2008), p.11-19.
  \bibitem[3]{popi}Imran Anwar, Dorin Popescu "Stanley Conjecture in small embedding dimension", Journal of Algebra 318 (2007),                            p.1027-1031.
  \bibitem[4]{apel}J.Apel "On a conjecture of R.P.Stanley", Journal of Algebraic Combinatorics, 17(2003), p.36-59.
  \bibitem[5]{par}Csaba Biro, David M.Howard, Mitchel T.Keller, William T.Trotter, Stephen J.Young "Partitioning subset lattices into intervals, preliminary version", Preprint 2008.   
  \bibitem[6]{mir}Mircea Cimpoeas "Stanley depth for monomial complete intersection", Bull. Math.
                  Soc. Sc. Math. Roumanie 51(99), no.3 (2008), p.205-211.%, http://arxiv.org/pdf/0805.2306.
  \bibitem[7]{mir2}Mircea Cimpoeas "Stanley depth for monomial ideals in three variables", Preprint 2008, http://arxiv.org/pdf/0807.2166.
  \bibitem[8]{hsy}J\"urgen Herzog, Ali Soleyman Jahan, Siamak Yassemi "Stanley decompositions and partitionable simplicial complexes",                    Journal of Algebraic Combinatorics 27(2008), p.113-125.
  \bibitem[9]{hvz}J\"urgen Herzog, Marius Vladoiu, Xinxian Zheng "How to compute the Stanley depth of a monomial ideal", to appear in                     Journal of Algebra%, http://arxiv.org/pdf/0712.2308.   
  \bibitem[10]{jah}Ali Soleyman Jahan "Prime filtrations of monomial ideals and polarizations", Journal of Algebra  312 (2007),                            p.1011-1032.
  \bibitem[11]{sum}Sumiya Nasir "Stanley decompositions and localization", 
                  Bull. Math. Soc. Sc. Math. Roumanie 51(99), no.2 (2008), p.151-158.  
  \bibitem[12]{pop}Dorin Popescu "Stanley depth of multigraded modules", Preprint 2008, http://arxiv.org/pdf/0801.2632.
  \bibitem[13]{asia}Asia Rauf "Stanley decompositions, pretty clean filtrations and reductions modulo regular elements", Bull. Math.
                    Soc. Sc. Math. Roumanie 50(98), no.4 (2007), p.347-354.
  \bibitem[14]{rin}Giancarlo Rinaldo "An algorithm to compute the Stanley depth of monomial ideals", Preprint 2008.
  \bibitem[15]{shen}Yihuang Shen "Stanley depth of complete intersection monomial ideals and upper-discrete partitions", Preprint 2008
  %,http://arxiv.org/pdf/0805.4461
\end{thebibliography}
\end{document}